\documentclass[12pt]{amsart}
\usepackage{amscd,amssymb}
\usepackage[arrow,matrix]{xy}

\theoremstyle{plain}

\newtheorem{prop}[subsection]{Proposition}
\newtheorem{cor}[subsection]{Corollary}

\theoremstyle{definition}
\newtheorem{rk}[subsection]{Remark}
\newtheorem{definition}[subsection]{Definition}
\newtheorem{ex}[subsection]{Example}

\numberwithin{equation}{section} 

\newcommand{\A}{{\mathcal A}}

\newcommand{\CC}{{\mathcal C}}

\newcommand{\Z}{\mathbb{Z}}

\newcommand{\C}{\mathbb{C}}

\newcommand{\PP}{\mathbb{P}}

\begin{document}

\title [On the torsion of Brieskorn modules of homogeneous polynomials]
{On the torsion of Brieskorn modules of homogeneous polynomials}

\author[ Shabbir Khurram]{ Shabbir Khurram }
\address{  68-B,New Muslim Town, School Of Mathematical Sciences GCU Lahore.
                 Pakistan}
 \email {khurramsms@gmail.com}
 \subjclass[2000]{Primary 14J70, 32S40 ; Secondary 14F45, 32S25.}

\keywords{torsion, Brieskorn module, Milnor algebra, homogeneous polynomial}

\begin{abstract}
Let $f\in \mathbb{C}[X_1,..., X_n]$ be a homogeneous polynomial
and $B(f)$ be the corresponding Brieskorn module. We describe the
torsion of the Brieskorn module $B(f)$ for $n=2$ and show that any torsion element has order 1. For $n>2$, we find some
examples in which the torsion order is strictly greater than 1.
\end{abstract}

\maketitle

\section {The Milnor algebra and the Brieskorn module}

Let $f\in R = \C[x_1,...,x_n]$ be a homogeneous polynomial of degree $d >1$. Then the
Koszul complex of the partial derivatives $f_j=\frac{\partial
f}{\partial x_j}$; $j=1,...,n$ in $R$ can be identified to the
complex
\begin{equation}\label{eq1}
0 \longrightarrow
\Omega^{0}\stackrel{df\wedge}{\longrightarrow}\Omega^{1}\stackrel{df\wedge}{\longrightarrow}
...\stackrel{df\wedge}{\longrightarrow}
\Omega^{n-1}\stackrel{df\wedge}{\longrightarrow}\Omega^{n}\longrightarrow
0
\end{equation}
where $\Omega^j$ denotes the regular differential forms of degree $j$ on
$\C^n$.

Let $J_f$ be the Jacobian ideal spanned by the partial derivatives $f_j$, $j=1,...,n$, in $R$
and $M(f)=R/J_f$ be the Milnor algebra of $f$.
One has the following obvious isomorphism of graded vector spaces
\begin{equation}\label{eq2}
M(f)(-n)= \frac{ \Omega^n}{df \wedge
\Omega^{n-1}}.
\end{equation}
Here, for any graded $\mathbb{C}[t]$-module $M$, the shifted module
 $M(m)$ is defined by setting $M(m)_s=M_{m+s}$ for all
 $s\in\mathbb{Z}.$

 We define the (algebraic) Brieskorn module as the quotient
\begin{equation}\label{eq3}
B(f)=\frac{ \Omega^n}{df \wedge d (\Omega^{n-2})}
\end{equation}
in analogy with the (analytic) local situation considered in \cite{Br}, see also \cite{L}, as well as
the submodule
\begin{equation}\label{eq4}
C(f) = \frac{df \wedge
\Omega^{n-1}}{df \wedge d (\Omega^{n-2})}.
\end{equation}
These modules are modules over the ring $\C[t]$ and the multiplication by $t$ is given by multiplying
by the polynomial $f$. Sometimes $B(f)$ is denoted by $G_f^{(0)}$ and $C(f)$ by $G_f^{(-1)}$, see \cite{BD},
\cite{DS1}.

 One has the following basic relation between the Milnor algebra and the Brieskorn
module, see \cite{DS2}, Prop. 1.6.

\begin{prop} \label{p1}
$$d f \wedge \Omega^{n-1} = df \wedge d(\Omega^{n-2})+f \cdot \Omega^n$$
In particular
$$B(f)/f.B(f) \simeq M(f)(-n)$$
\end{prop}
Let $B(f)_{tors}$ be the submodule of $\C[t]$-torsion elements in $B(f)$ and define the reduced Brieskorn module
${\overline B}(f)=B(f)/B(f)_{tors}$.

\begin{rk}\label{r0}
The reduced Brieskorn module
${\overline B}(f)$ is known to be a free $\C[t]$-module of rank $b_{n-1}(F)$, where $F=\{x \in \C^n ~~|~~f(x)=1\}$ is the affine Milnor fiber of $f$. Indeed, it follows from the section (1.8) in \cite{DS2}
that one has a canonical isomorphism $\oplus_{j=0,d-1}{\overline B}(f)_{qd+j} = H^{n-1}(F,\C)$
for any $q\ge n$. If we assume the $\C[t]$-basis of ${\overline B}(f)$ to be formed by homogeneous elements
(which is always possible), each basis element will contribute by 1 to the dimension of
$\oplus_{j=0,d-1}{\overline B}(f)_{qd+j}$ for $q$ large enough.

 Moreover, it exists an integer $N>0$ such that $t^N
B(f)_{tors}=0$, see \cite{DS2}, Remark 1.7. The least $N$ satisfying this condition is called the
{\it torsion order} of $f$ and is denoted by $N(f)$.
\end{rk}

\begin{rk}\label{r01} As noted in \cite{DS2}, Remark 1.7, there is a slight difference between the Brieskorn module defined above and the Brieskorn module considered by Barlet and Saito in \cite{BaS}.
In fact, the latter one is defined to be the $n$-th cohomology group $H^n(\A_f^*)$, where
$\A_f^j=\ker (df \wedge :\Omega ^j \to \Omega ^{j+1})$ and the differential $d_f:\A_f^j \to \A_f^{j+1}$
is induced by the exterior differential $d:\Omega ^j \to \Omega ^{j+1}$. Since
$ df \wedge d(\Omega^{n-2})\subset  d_f(\A_f^{n-1}) $, one has an epimorphism
$$B(f) \to H^n(\A_f^*),$$
(note that the direction of this arrow is misstated in \cite{DS2}, Remark 1.7).
Moreover, when $n=2$, it follows from Proposition 3.7 in \cite{BaS} that $H^2(\A_f^*)$ is torsion-free.
Our results in section 2 show that the torsion module $B(f)_{tors}$ can even be not of finite type in this setting, and hence
the above epimorphism is far away from an isomorphism in general.

\end{rk}

\begin{definition} \label{d1}

 For $b\in B(f)$, we say that $b$ is $t$-torsion of
order $k_b \geq 1$ if $t^{k_b}\cdot b=0$ and $t^{k_b-1}\cdot b\neq0.$
\end{definition}

It is clear that such a torsion order $k_b$ divides the torsion order $N(f)$ of $f$, since $t^{N(f)}\cdot b=0$
and that $N(f)$ is the G.C.D. of all the $t$-torsion orders $k_b$ for $b \in B(f)$ a torsion element.

\subsection {The case of an isolated singularity}

Assume that $f\in \mathbb{C}[x_1,...,x_n]$ is a homogeneous
polynomial having an isolated singularity at the origin of $\C^n$.
Then the dimension of the Milnor algebra $M(f)$ (as a $\C$-vector
space) is the Milnor number of $f$ at the origin, denoted by
$\mu(f)$. One has in this case $\mu(f)=b_{n-1}(F)$, see for instance
\cite{D1}. In this case, the structure of the module $B(f)$ is as
simple as possible.

 \begin{prop} \label{p2}
  The $\C[t]$-module $B(f)$ is free of rank
 $\mu(f).$
 \end{prop}

 \proof

 Indeed, a homogeneous polynomial $f$ having an
 isolated singularity at the origin induces a tame mapping $f:\mathbb{C}^n\rightarrow \mathbb{C}$,
 to which the results in \cite{DS1}, \cite{NS} and \cite{S} apply.
 Indeed, the Gauss-Manin system $G_f$ of $f$, which is an $A_1=\C[t]<\partial_t>$-module, contains a
 $\C[t]$-submodule $G_f^{(0)}$, which is known to be free of finite type for a tame polynomial, see
 for instance Remark 3.3 in \cite{DS1}. As we have already mentionned above, $G_f^{(0)}=B(f)$.

 \endproof

 \begin{cor} \label{c1}
There is an isomorphism
$$B(f)(n)=M(f)\otimes_\mathbb{C}\mathbb{C}[t]
$$
of graded modules over the graded ring
$\mathbb{C}[t]$.
 In particular, one has, at the level of the associated Poincar\'e series, the following equality.
$$P_{B(f)}(t)=P_{M(f)}(t)\cdot \frac{t^n}{1-t}=\frac{t^n\cdot (1-t^{d-1})^n}{(1-t)^{n+1}}.$$
\end{cor}

\begin{ex} \label{e1}

For $n=3,$ take $f(x,y,z)=x^3+y^3+z^3.$ Then a
$\mathbb{C}$-basis of $M(f)$ is given by $1,x,y,z,xy,yz,xz,xyz.$ The
same $8$ monomials form a basis of $B(f)$ as a free
$\mathbb{C}[t]$-module.
In this case
$$P_{B(f)}(t)=\frac{t^3(1-t^{2})^3}{(1-t)^{4}}=t^3+4t^4+7t^5+8t^6+8t^7+8t^8+\cdot\cdot\cdot.$$

\end{ex}

 \begin{cor} \label{c2}
 The $\mathbb{C}[t]$-module $B(f)$ is torsion free
 if and only if $0$ is an isolated singularity of the homogeneous polynomial
 $f$.
 \end{cor}

 \proof If $0$ is not an isolated singularity, then the Milnor algebra $M(f)$ is an infinite dimensional
$\C$-vector space. The isomorphism in Proposition \ref{p1} implies that in this case $B(f)$ is not
finitely generated over $\C[t]$. By Remark \ref{r0}, it follows that the canonical projection
$$B(f) \to {\overline B}(f)$$
is not an isomorphism, hence $B(f)_{tors} \ne 0.$

\endproof

Proposition \ref{p1} implies that $f \cdot B(f)=C(f)$, which in turn yields the following.

\begin{cor} \label{c2.5}

 Assume that $0$ is not an isolated singularity of the homogeneous polynomial
 $f$. Then $N(f)=1$ if and only if the $\mathbb{C}[t]$-module $C(f)$ is torsion free.
 \end{cor}

 \proof
  Let $b\in B(f)_{tors}$ be a non-zero element, which exists by our assumption. By  Remark \ref{r0}, it follows that $b$ is $t$-torsion, say of order
 $k.$ If $k>1,$ then $$0=t^k\cdot b=t^{k-1}\cdot(tb).$$
 By Proposition \ref{p1},  we know that $t\cdot b \in C(f).$ If $C(f)$ is torsion free,
 we get that $t\cdot b=0$ in $C(f)$, i.e. $t\cdot
 b=0$ in $B(f)$, a contradiction. Hence $k=1$ for any $b\in B(f)_{tors}$, in other words $N(f)=1$.

 \endproof

 \section {The case  $n=2$}
 In this section we suppose that $f\in\mathbb{C}[x,y]$ is a
 homogeneous polynomial of degree $d>1$, which is not the power $g^r$
 of some other polynomial $g\in\mathbb{C}[x,y]$ for some $r>1.$ This
 condition is equivalent to asking the Milnor fiber $F$ of $f$ to be
 connected and such polynomials are sometimes called {\it primitive}. For more on this, see \cite{DP}, final Remark, part (I).

\begin{prop} \label{p3}

 The submodule $C(f)$ is a free $\C[t]$-module of rank $b_1(F)$.
\end{prop}

\proof

The fact the $C(f)$ is torsion free follows from Proposition 7. $(ii)$ in
\cite{BD}. Indeed, it is shown there that the localization morphism
$$\ell:C(f) \to C(f) \otimes_{\C[t]}\C[t,t^{-1}]$$
is injective, which is clearly equivalent to the fact that there are no $t$-torsion elements.

 For the second claim, note that the composition of the inclusion $C(f) \to B(f)$ and the
canonical projection
$B(f) \to {\overline B}(f)$ gives rise to an embedding of $C(f)$ into ${\overline B}(f)$ whose image is exactly $t \cdot {\overline B}(f)$ (use again Proposition \ref{p1} and the fact that $C(f)$ is torsion free). Since $\C[t]$ is a principal ideal domain, the claim follows from the structure theorem of submodules  of free modules
of finite rank over such rings.

\endproof

Corollary  \ref{c2.5} implies the following.

 \begin{cor} \label{c3}
 If $f \in \C[x,y]$ is a homogeneous polynomial with a non-isolated singularity at the origin,
 then $N(f)=1$.
  \end{cor}

 The above Corollary can be restated by saying that
 $$0\rightarrow B(f)_{tors}\rightarrow B(f)\stackrel{t}\rightarrow C(f)\rightarrow 0$$
 is
 an exact sequence of graded $\C[t]$-modules.
 We get thus an
 isomorphism of graded $\mathbb{C}[t]$-modules
 $$\overline{B}(f)(-d)\simeq C(f).$$

 \begin{ex} \label{e2}
 Let $f=x^py^q$ with $(p,q)=1$. In order to compute the torsion of the Brieskorn module $B(f),$ we start by
  finding  $\C$-vector space monomial bases for $\frac{B(f)}{C(f)}\simeq M(f)$ (up-to a shift in gradings) and $C(f)$. Note that the Jacobian ideal is given by
  $$J_f=\langle x^{p-1}y^q,x^py^{q-1}\rangle=x^{p-1}\cdot y^{q-1}\langle x,y\rangle.$$
It follows that $\frac{B(f)}{C(f)}\simeq \frac{S}{J_f}$ is an
infinite dimensional $\mathbb{C}$-vector space with a monomial basis
given by $x^ay^b$ with $a\leq p-2$ or $b\leq q-2$ or
$(a,b)=(p-1,q-1).$

We have to compute $df \wedge d\Omega^0=\{df \wedge dg\} ,$ where $g$ is a polynomial
function on $\mathbb{C}^2$.
Since we are working with homogeneous polynomials, it is enough to work with one monomial, say $(x^ay^b),$ at a time. We have
$$df\wedge d(x^{a}y^{b})=
(pb-qa)x^{a+p-1}y^{b+q-1}dx\wedge dy.$$
Since $$C(f)=\frac{J_f\cdot\Omega^2}{df\wedge d\Omega^0
}=\frac{x^{p-1}y^{q-1}\langle x,y\rangle \Omega^2}{df\wedge
d\Omega^0},$$
 a system of generators of the
$\mathbb{C}$-vector space $C(f)$ is given by the
classes of the elements $w\in x^{p-1}y^{q-1}\langle x,y\rangle
\Omega^2,$ which do not belong to $df\wedge d\Omega^0.$
 To find  them, it is enough to look at monomial differential forms i.e. $x^\alpha y^\beta dx \wedge dy$ with $\alpha\geq p -1$, $ \beta\geq q -1$ and $\alpha +\beta \ge p+q-1.$

 So,
if $a+p-1=\alpha, b+q-1=\beta$ and $pb-qa \neq 0 ,$  then
$$df \wedge d(\frac{x^a y^b}{pb-aq})= x^\alpha y^\beta dx \wedge dy.$$
Hence, the only elements  $x^\alpha y^\beta dx \wedge dy$ which are
not in $df \wedge d \Omega ^0$
 are $x^{a+p-1}y^{b+q-1}dx \wedge dy$ where $pb=qa$, i.e. $ a=kp,b=kq$ for some $ k\geq 0.$
It follows that $C(f)$, as a $\mathbb{C}$-vector space, has a basis
given by $x^{(k+1)p-1}y^{(k+1)q-1}$ where $k\geq 0$,  which can be
written as
$$C(f)=\C[t]\cdot x^{2p-1}y^{2q-1} dx \wedge dy$$
i.e. $C(f)$ is a free $\C[t]$-module of rank 1.

Moreover, the corresponding monomial basis as a $\C$-vector space for the Brieskorn module $B(f)$ is given by $x^a y^b dx \wedge dy$ with
$ a\leq p-2$ or $ b\leq q-2$ or $a={(k+1)p-1}, b={(k+1)q-1}$ for some $ k\geq 0.$

With respect to this basis,  $B(f)_{tors}$ is the linear span of $ x^ay^b dx \wedge dy$ with
$a\leq p-2$ or $b\leq q-2$. Indeed our computation
above shows that
$$df\wedge d(x^{a+1}y^{b+1})=(p(b+1)-q(a+1))f\cdot w$$
 i.e. $tw=0$ if the coefficient $p(b+1)-q(a+1)$ is non-zero.

 It follows that
$$\overline{B}(f)=\mathbb{C}[t]\cdot [x^{p-1}y^{q-1}dx\wedge dy],$$
a
free $\mathbb{C}[t]$-module of rank 1.

It remains to explain why  $b_1(F)=1$, where $F=\{ (x,y)\in
\mathbb{C}^2 ; x^py^q=1\}.$
Consider the covering
$$\mathbb{Z}_q \hookrightarrow F \stackrel{\phi}\rightarrow \mathbb{C}^\ast$$
where $\phi $ maps $(x,y)$ to $x.$\\
Now a covering yields an exact sequence (see for
instance \cite{AH},page 376),
$$0 \hookrightarrow
 \pi_1(F,(x_0,y_0))\stackrel{\phi_\#} \rightarrow \pi_1
(\mathbb{C}^\ast,x_0)\rightarrow
\pi_0(\mathbb{Z}_q,(x_0,y_0))=\Z_q \to 0.$$
 Hence $\pi_1(F)\simeq \mathbb{Z}$, which
shows that the first Betti number $b_1(F)$ is $1.$
 \end{ex}

\section {Eigenvalues of the monodromy and torsion of Brieskorn modules}

Any homogeneous polynomial $f \in \C[x_1,...,x_n]$ induces a locally trivial fibration
$f:\C^n \setminus f^{-1}(0) \to \C^*$, with fiber $F$ and semisimple monodromy operators
$$T_f^k:H^k(F,\C) \to H^k(F,\C)$$
for $k=0,...,n-1.$
The eigenvalues of $T_f$ are $d$-roots of unity, where $d$ is the degree of $f$,
and for each such eigenvalue $\lambda$ we denote by $H^*(F,\C)_\lambda$
the corresponding eigenspace.

According to \cite{DS2}, see the discussion just before Remark 1.9, one has for
$q < n$ an inclusion
$$t^{n-q}: \overline{B}(f)_{qd-j}\rightarrow
\overline{B}(f)_{nd-j}$$
and, for $q \ge n$, an identification
$$\overline{B}(f)_{qd-j}  =H^{n-1}(F,\mathbb{C})_\lambda $$
where $\lambda = \exp(\frac{2\pi j \sqrt{-1}}{d}),$ with
$j=0,1,...,d-1.$
Let $\omega_n=dx_1\wedge...\wedge dx_n$ and note that $[\omega _n] \in \overline{B}(f)_{qd-j}$ for some $q$
if and only if $n-j$ is divisible by $d$.
This yields the following, see \cite{DS2}, Corollary 1.10.

\begin{cor} \label{c4}
 Assume that $\lambda =\exp(\frac{2\pi n \sqrt{-1}}{d})$ is
 not an eigenvalue of the monodromy operator $T_f$ acting on $H^{n-1}(F,\C)$. Then  $[\omega_n]$ is a non-zero
 torsion element
 in $B(f).$
\end{cor}

 Here are some examples of torsion elements in $B(f)$ obtained using this approach, and having torsion order $N(f)>1$.

\begin{ex} \label{e3}

 Let $f=x^3+y^2z,$  $n=3,\,\,d=3,\,\,j=0 $ (this defines
 a cuspidal cubic
curve $\CC$ in $\PP^2$). In this case $\lambda =1$ is not an eigenvalue of the monodromy acting on
$H^2(F,\C)$. Indeed, one has $H^2(F,\C)_1=H^2(U,\C)=0$. Here $U =\PP^2 \setminus \CC$ and the last vanishing comes from the following obvious equalities:
$b_0(U)=1$, $b_1(U)=0$ and $\chi(U)=\chi (\PP^2)- \chi (\CC)=3-2=1$ as $\CC$ is homeomorphic to a sphere $S^2$.

For $w=dx\wedge  dy\wedge dz,$ to find the value of $k$ such that
$t^k\cdot [w]=0$ in $B(f)$ is
equivalent to saying that $f^k\cdot w \in df \wedge d\Omega^1$.
We have to check  for which value of $k,$ we have
solutions of the equation:
$$f^k\cdot w =[3x^2(\frac{\partial R}{\partial y}-\frac{\partial Q}{\partial z})+
2yz(\frac{\partial P}{\partial z}-\frac{\partial R}{\partial x})+
y^2 (\frac{\partial Q}{\partial x}-\frac{\partial P}{\partial
y})]dx\wedge  dy\wedge dz$$ where $P,Q,R\in \mathbb{C}[x,y,z]$ are homogeneous polynomials of degree $3k-1$.

For $k=1,$ we  get a system of non-homogeneous linear equations in
which we have $15$ unknowns and $9$ equations. Then using  the software Mat
Lab we compute the rank of the matrix corresponding to the homogeneous system (containing
$9$ rows and $15$ columns) and get  $8.$ On other hand, the rank of the matrix associated to the
non-homogeneous system (containing $9$ rows and $16$ columns) is
$9,$ which show
that this system has no solution.

For $k=2,$  we  get another system of non-homogeneous linear equations
containing $60$ unknowns and $26$ equations. Then using Mat Lab we
compute the rank of the matrix corresponding to the homogeneous system (containing $26$
rows and $60$ columns) and get  $26.$ This time the corresponding
non-homogeneous system matrix (containing $26$ rows and $61$ columns) has also
rank $26,$ which shows that this system of equations has a solution.
An explicit solution for $k=2,$ is $P=\frac{32}{3}( x^3yz),$
$Q=xy^2z^2$ and $R=\frac{1}{3}( x^4y).$
Hence  $[w]$ has $t$-torsion order 2 in $B(f).$

 \end{ex}

 The next example shows somehow what happens when $H^{n-1}(F,\C)_\lambda \ne 0 $.

\begin{ex} \label{e4}
 Let $f=x^2y^2+xz^3+yz^3$ where  $n=3,\,\,d=4 $ and $j=3$. The corresponding curve $\CC:f=0$ in $\PP^2$
 has  two cusps as singularities. It follows from the study of the plane quartic curves, see \cite{D1}, p. 130, that $\pi_1(U)=\Z_4$ and hence $H^1(F,\C)=0$.

 Next $\chi(U)=\chi(\PP^2)- \chi(\CC)=3-(2-2\cdot 3+2\cdot 2)=3$. The zeta-function of the monodromy operator $T_f$ looks like
 $$\det(t\cdot Id-T_f^0)\cdot \det(t\cdot Id-T_f^2)=(t^4-1)^3$$
 see for instance \cite{D1}, p. 108.
 Hence we can not apply Corollary \ref{c4} to this polynomial to infer that $[w]=[\omega_3]$ is a torsion element in $B(f)$.

However, we can try to find values of $k,$ for which we have
solutions of the equation:
$$f^k\cdot w =[(2xy^2+z^3)(\frac{\partial R}{\partial y}-\frac{\partial Q}{\partial z})+
(2x^2y+z^3)(\frac{\partial P}{\partial z}-\frac{\partial R}{\partial
x})+ (3xz^2+3yz^2) (\frac{\partial Q}{\partial x}-\frac{\partial
P}{\partial
y})]w$$
 where $P,Q,R\in \mathbb{C}[x,y,z]$ are homogeneous polynomials of degree $4k-2$.



For $k=1$ and $k=2$ a similar computation of matrix ranks as above shows that the corresponding systems have no solution.
Hence  $[w]$ is either a non-torsion element, or it has $t$-torsion order greater or equal to $3$ in $B(f).$

 \end{ex}

 \begin{rk} \label{r1}
 Note that in both examples above the element $t\cdot [w]$ is a non-zero torsion element in $C(f)$.
 Hence Proposition \ref{p3} fails for $n>2$. It also shows that in these cases the module $C(f)$
 is not torsion free, compare to Corollary \ref{c2.5}
\end{rk}

The last example shows that even for rather complicated examples (here the zero set of $f$ is a surface $S$ with non-isolated singularities) one may still have 1 as the $t$-torsion order of $[\omega_n]$.

\begin{ex} \label{e5} Let $f=x^2z+y^3+xyt,$ $n=4,\,\,d=3,\,\,j=1$ be the
equation of a cubic surface $S$  in $\PP^3$. It follows from \cite{D0}, Example 4.3, that  $H^3(F,\C)=0$.
Hence we can  apply Corollary \ref{c4} to this polynomial to infer that $[w]=[\omega_4]$ is a torsion element in $B(f)$.

We have to check  for which values of $k,$ we have
solutions of the equation:
$$f^k\cdot w =[(2xz+yt)(\frac{\partial S}{\partial t}-\frac{\partial
T}{\partial z} +\frac{\partial U}{\partial
y})+(3y^2+xt)(-\frac{\partial Q}{\partial t}+\frac{\partial
R}{\partial z} -\frac{\partial U}{\partial x})+$$
$$+ x^2 (\frac{\partial
P}{\partial t}-\frac{\partial R}{\partial y} +\frac{\partial
T}{\partial x})+ xy (\frac{\partial P}{\partial z}-\frac{\partial
Q}{\partial y} -\frac{\partial S}{\partial x})]dx\wedge dy \wedge dz
\wedge dt,$$
 where $P,Q,R,S,T,U \in \mathbb{C}[x,y,z,t]$ are homogeneous polynomials of degree $3k-1$.

For $k=1,$ we  get a system of non-homogeneous linear equations in
which we have $42$ unknowns and $15$ equations. Using MatLab we compute the rank of the corresponding homogeneous system (containing
$15$ rows and $42$ columns) and get $14.$ And the rank of non-homogeneous
system (containing $15$ rows and $43$ columns) has the same rank $14,$
which
show that this system has a solution.

Hence  $t\cdot [w]=0$ in $B(f).$

\end{ex}

\end {document}